\documentclass{amsart}

\usepackage{amssymb,amscd,amsthm,verbatim,psfig,epsfig
}

\newtheorem{proposition}{Proposition}
\newtheorem{theorem}[proposition]{Theorem}

\newtheorem{lemma}[proposition]{Lemma}

\theoremstyle{remark}

\def\ignore#1{\relax}

\title[Quantum Teichm\"{u}ller
Spaces] {A Uniqueness Property for the Quantization\\ of Teichm\"{u}ller
Spaces}
\author{Hua Bai}
\email{huabai@usc.edu}
\address{Department of Mathematics\\
           University of Southern California\\
           Los Angeles, CA~90089-2532, U.S.A.}


\thanks{ This work was partially
supported by the grant
DMS-0103511 from the National
Science Foundation}
\date{\today}

\begin{document}

\begin{abstract}
Chekhov, Fock and Kashaev introduced  a quantization of the
Teichm\"{u}ller space $\mathcal{T}^q(S)$ of a punctured surface $S$, and
an exponential version of this construction was developed by
Bonahon and Liu. The construction of the quantum Teichm\"{u}ller
space crucially depends on certain coordinate change isomorphisms
between the Chekhov-Fock algebras associated to different ideal
triangulations of $S$. We show that these
coordinate change isomorphisms are essentially unique, once we require
them to satisfy a certain number of natural conditions.
\end{abstract}

\maketitle

Let $S$ be an oriented surface of finite topological type, with at
least one puncture. The \textit{Teichm\"{u}ller space}
$\mathcal{T}(S)$ is the space of isotopy classes of complete
hyperbolic metrics on $S$. A quantization of the Teichm\"{u}ller
space $\mathcal{T}(S)$ of $S$ was introduced by L.~Chekhov and V.~Fock
\cite{ChekhovFock99,ChekhovFock00,Fock97} and, independently, by
R.~Kashaev \cite{Kashaev98} (see also \cite{Teschner04}) as an approach to
quantum gravity in $2+1$ dimensions. This is a deformation of the
$\mathrm C^*$-algebra of functions on the usual Teichm\"{u}ller space
$\mathcal{T}(S)$ of $S$, depending on a parameter $\hbar$, in such a way
that the linearization of this deformation at $\hbar = 0$ corresponds to
the Weil-Petersson Poisson structure of $\mathcal{T}(S)$. F.~Bonahon and
X.~Liu \cite{BonahonLiu04,Liu04} developed an exponential version of the
Chekhov-Fock-Kashaev construction. This exponential
version of the quantization can be formulated in terms
of non-commutative algebraic geometry, and has the advantage
of possessing an interesting finite dimensional representation
theory \cite{BonahonLiu04}, whereas the non-exponential version is
defined in terms of self-adjoint operators of Hilbert spaces.

More precisely, let $S$ be an oriented punctured surface of finite
topological
   type, obtained by removing a finite set $\{v_1,v_2,\dots,v_p\}$
   from a closed oriented surface $\overline{S}$.
   An \textit{ideal triangulation} is a family $\lambda$ of finitely
many disjoint simple arcs $\lambda_1,\lambda_2,\dots,\lambda_n$
going from puncture to puncture and decomposing $S$ into finitely
many triangles with vertices at infinity; in other words, an ideal
triangulation consists of the edges of a triangulation of the
closed surface $\overline S$ whose vertex set consists of the
punctures $\{v_1,v_2,\dots,v_p\}$. Considering $q=e^{\pi i \hbar}$
as an indeterminate over $\mathbb{C}$, the \textit{Chekhov-Fock
algebra} $\mathbb{C}[X_1,X_2,\dots,X_n]_{\lambda}^q$ associated to
the ideal triangulation $\lambda$ is the algebra over
$\mathbb{C}(q)$ defined by generators
$X_1^{\pm1},X_2^{\pm1},\dots,X_n^{\pm1}$ associated to the
components of $\lambda$, and by relations $X_i X_j X_i^{-1}
X_j^{-1} = q^{2\sigma_{ij}} $ where the $\sigma_{ij}$ are integers
determined by the combinatorics of the ideal triangulation and
connected to the Weil-Petersson form on Teichm\"uller space. This
algebra has a well-defined fraction division algebra
$\mathbb{C}(X_1,X_2,\dots,X_n)_{\lambda}^q$. In practice, the
algebra $\mathbb{C}(X_1,X_2,\dots,X_n)_{\lambda}^q$ consists of
all formal rational fractions in variables $X_i$ that skew-commute
according to the relations $X_i X_j = q^{2\sigma_{ij}} X_j X_i$.

   As one moves from one ideal triangulation $\lambda$  to another
$\lambda'$, Chekhov and
Fock \cite{Fock97, ChekhovFock00, ChekhovFock99} (as developed in
\cite{Liu04, BonahonLiu04}) introduce explicit \textit{coordinate change
isomorphisms}
$$\Phi_{\lambda
\lambda'}^q:\mathbb{C}(X_1,X_2,\dots,X_n)_{\lambda'}^q
\rightarrow \mathbb{C}(X_1,X_2,\dots,X_n)_{\lambda}^q.
$$
These are algebra isomorphisms which satisfy the natural property
that $\Phi_{\lambda'' \lambda'}^q \circ \Phi_{\lambda' \lambda}^q
=\Phi_{\lambda'' \lambda}^q $ for any ideal triangulations
$\lambda$, $\lambda'$, $\lambda''$. In a triangulation independent
way, this associates to the surface $S$ the algebra
$\mathcal{T}^q(S)$ defined as the quotient of the family of all
$\mathbb{C}(X_1,X_2,\dots,X_n)_{\lambda}^q$, with $\lambda$
ranging over all the ideal triangulations of the surface $S$, by
the equivalence relation that identifies
$\mathbb{C}(X_1,X_2,\dots,X_n)_{\lambda'}^q$ and
$\mathbb{C}(X_1,X_2,\dots,X_n)_{\lambda}^q$
   by the coordinate change isomorphism $\Phi_{\lambda
\lambda'}^q$ . By definition,  $\mathcal{T}^q(S)$ is the
\textit{quantum Teichm\"{u}ller space} of the surface $S$.

This construction is motivated by the case where $\hbar =0$, so that $q =
\mathrm e^{2\pi i \hbar}=1$. Thurston associated to each ideal
triangulation $\lambda$ a set of global coordinates for the (non-quantum)
Teichm\"uller space $\mathcal{T}(S)$, called shear coordinates; see
\cite{Thu, Fock97, Bonahon97}. When $q=1$,
$\mathbb{C}(X_1,X_2,\dots,X_n)_{\lambda}^1$ is the usual algebra
$\mathbb{C}(X_1,X_2,\dots,X_n)$ of rational functions in the commuting
variables $X_i$, and  the coordinate change isomorphism
$\Phi_{\lambda
\lambda'}^1:\mathbb{C}(X_1,X_2,\dots,X_n)
\rightarrow \mathbb{C}(X_1,X_2,\dots,X_n)
$
exactly corresponds to the coordinate change between the shear
coordinates $X_i$ for $\mathcal{T}(S)$ respectively associated to
$\lambda$ and $\lambda'$. As a consequence,
$\mathcal{T}^1(S)$ is in a natural way the
algebra of rational functions on the  Teichm\"uller space
$\mathcal{T}(S)$.

We need to be a little more specific in our definitions, by requiring the
data of an ideal triangulation
$\lambda$ to include an indexing of its components
$\lambda_1$, $\lambda_2$, \dots, $\lambda_n$ by the set
$\{1,2,\dots,n\}$. The permutation group $\mathcal{S}_n$ then acts by
reindexing on the set $\Lambda(S)$ of isotopy classes of all such indexed
ideal triangulations of $S$.

The coordinate change isomorphisms
defined by Chekhov-Fock \cite{ChekhovFock99}, Kashaev
\cite{Kashaev98} and Liu
\cite{Liu04} satisfy the following natural conditions:

   \begin{theorem}[Chekhov-Fock, Kashaev, Liu]
\label{thm:CoordChangesExist}
   There exists a family of algebra isomorphisms
$$
\Phi_{\lambda
\lambda'}^q:\mathbb{C}(X_1',X_2',\dots,X_n')_{\lambda'}^q
\rightarrow \mathbb{C}(X_1,X_2,\dots,X_n)_{\lambda}^q
$$
indexed by pairs of ideal triangulations $\lambda$, $\lambda'\in
\Lambda(S)$, which satisfy the following conditions:
\begin{enumerate}
   \item[(1)]
$\Phi_{\lambda \lambda''}^q=\Phi_{\lambda \lambda'}^q \circ
\Phi_{\lambda' \lambda''}^q$  for any $\lambda,\lambda'$ and
$\lambda'' \in \Lambda(S)$ ;
   \item[(2)]if $\lambda'=\sigma
\lambda$ is obtained by reindexing $\lambda$ by a permutation $\sigma \in
\mathcal S_n$, then
$\Phi_{\lambda \lambda'}^q(X_i')=X_{\sigma(i)}$  for any
$1\leqslant i\leqslant n$;
   \item[(3)] a Locality Condition precisely described in
Section~\ref{sect:LocalityCondition}.
\end{enumerate}
   \end{theorem}

This paper is devoted to a uniqueness property for these  $\Phi_{\lambda
\lambda'}^q$. This will require the property that the surface
$S$ is sufficiently large in the sense that its Euler
characteristic $\chi(S)$ is less than
$-2$. This excludes the sphere with $\leq 4$ punctures, and the torus
with $\leq 2$ punctures.

\begin{theorem}
\label{thm:CoordChangesUnique}
   Assume that the surface $S$ is sufficiently large, in the sense that
$\chi(S)<-2$. Then the family of coordinate change isomorphisms
$\Phi_{\lambda\lambda'}^q$
in Theorem~\ref{thm:CoordChangesExist} is unique up to a uniform rescaling
and/or inversion of the $X_i$.

Namely, if
$$
\Psi_{\lambda
\lambda'}^q:\mathbb{C}(X_1',X_2',\dots,X_n')_{\lambda'}^q
\rightarrow \mathbb{C}(X_1,X_2,\dots,X_n)_{\lambda}^q
$$ is another family of
isomorphisms satisfying the conditions of
Theorem~\ref{thm:CoordChangesExist}, then there exists a non-zero
constant $\xi \in \mathbb{C}(q)$ and a sign $\varepsilon =\pm1$
such that $\Psi_{\lambda \lambda'}^q=\Theta_{\lambda} \circ
\Phi_{\lambda \lambda'}^q \circ \Theta_{\lambda'}^{-1}$ for any
pair of ideal triangulations $\lambda$, $\lambda'$, where
$\Theta_{\lambda}: \mathbb{C}(X_1,X_2,\dots,X_n)_{\lambda}^q
\rightarrow \mathbb{C}(X_1,X_2,\dots,X_n)_{\lambda}^q$ is the
isomorphism defined by $\Theta_{\lambda}(X_i)=\xi X_i^\varepsilon$
for every $i$.
\end{theorem}

Theorem~\ref{thm:CoordChangesUnique} is false when $S$ is the
once-punctured torus or the 3--times punctured sphere. The
uniqueness property appears to hold for the twice-punctured torus and
the 4--times punctured sphere.

The proof of Theorem~\ref{thm:CoordChangesUnique} relies on two
crucial ingredients. One is that the $\Psi_{\lambda
\lambda'}^q $ are algebra isomorphisms. The other one is that they
satisfy the pentagon relation discussed in Section~\ref{ProofMainStep}.

This uniqueness property has a
counterpart in the non-exponential context of \cite
{Kashaev98, ChekhovFock99, ChekhovFock00, Fock97, Teschner04}. In that
context, it is crucial that the coordinate change isomorphisms
are defined by transforming self-adjoint operators under meromorphic
functions satisfying appropriate conditions. The punch line is then
provided by the characterization of the meromorphic function
$$
\phi^\hbar (z) = - \frac{\pi\hbar}2 \int_{-\infty}^{\infty}
\frac{\mathrm e ^{- \mathrm i tz}}{\sinh \pi t \, \sinh \pi \hbar
t}\,dt
$$
by a certain functional equation. By contrast, our arguments are purely
algebraic.

Our proof strongly uses the property that the algebras considered
are defined over
$\mathbb C(q)$, or at least that $q$ is not a root of unity. On the
other hand, interesting topological
applications occur when we specialize $q$ to a root of unity
\cite{BonahonLiu04}. It would be interesting to investigate what kind
of uniqueness property holds when $q^N=1$, including
the classical case $q=1$.

\medskip
\noindent\textit{Acknowledgments:} The author would like to thank Francis
Bonahon and Xiaobo Liu for many conversations related to this work.

\section{The Chekhov-Fock algebra and coordinate change isomorphisms}
\label{sect:CheFock}
In the surface  $S=\overline S- \{v_1, v_2,
\dots, v_p\}$, let $\lambda$ be an ideal triangulation consisting of
finitely many disjoint simple arcs $\lambda_1,\lambda_2,\dots,\lambda_n$
going from puncture to puncture. An easy computation shows that the
number $n$ of arcs in $\lambda$ is equal to
$-3\chi(S)=6g+3p-6$, where $\chi(S)$ is the Euler characteristic
of $S$,  $g$ is its genus and $p$ is its number of punctures.

   The complement $S-\lambda$ has $2n$ spikes converging towards the
punctures, and each spike is delimited by one $\lambda_i$ on one side and
one $\lambda_j$ on the other side, with possibly $i=j$. Let $a_{ij}$
denote the number of spikes of $S-\lambda$ that are delimited on the
left by
$ \lambda_i$ and on the right by $\lambda_j$, and define
$\sigma_{ij}=a_{ij}-a_{ji}$. It is immediate that
$\sigma_{ij}\in \{0,\pm1,\pm2\}$ and $\sigma_{ij}=-\sigma_{ji}$.

Let $q$ be an  indeterminate over $\mathbb{C}$. The \textit{Chekhov-Fock
algebra} $\mathbb{C}[X_1, X_2, \dots, X_n]_{\lambda}^q$, as defined in
\cite{BonahonLiu04, Liu04}, is the algebra over
$\mathbb{C}(q)$ defined by generators
$X_1^{\pm1}$, $X_2^{\pm1}$, \dots, $X_n^{\pm1}$ associated to the
components of $\lambda$ and by skew-commutativity relations
$$
X_i X_j X_i^{-1}
X_j^{-1} = q^{2\sigma_{ij}}.
$$

This is an iterated  skew-polynomial algebra. As such, it satisfies the
Ore Condition and consequently admits a well-defined fraction
division algebra
$\mathbb{C}(X_1,X_2,\dots,X_n)_{\lambda}^q$; see for instance
\cite{BrownGoodearl02, Cohn95, Kassel95}. In concrete terms, the
Chekhov-Fock algebra
$\mathbb{C}[X_1,X_2,\dots,X_n]_{\lambda}^q$ consists of
formal Laurent polynomials in variables $X_i$ satisfying the
skew-commutativity relations $X_i X_j X_i^{-1} X_j^{-1} =
q^{2\sigma_{ij}} $, while its fraction algebra
$\mathbb{C}(X_1,X_2,\dots,X_n)_{\lambda}^q$ consists of all formal
rational fractions in the $X_i$ satisfying the same relations.

\section{The Locality Condition}
\label{sect:LocalityCondition}
   We describe here the Locality Condition mentioned in
   the introduction.

   Define the \textit{discrepancy span} $D(\lambda,
\lambda')$ of two ideal triangulations $\lambda$, $\lambda'$ as
the closure of the union of those connected components of
$S-\lambda$ which are not isotopic to a component of $S-\lambda'$.
For instance, $D(\lambda, \lambda')$ is empty exactly when
$\lambda$ and $\lambda'$ coincide after isotopy and reindexing.

   The coordinate change
isomorphisms $\Phi_{\lambda \lambda'}^q$ are said to satisfy the
\textit{Locality Condition} if the following holds. Let $\lambda$ and
$\lambda'$ be two ideal triangulations indexed in such a way that
$\lambda_i \subset D(\lambda,
\lambda')$ when  $i \leq k$, and $\lambda_i' = \lambda_i$ when
$i>k$. Then the Locality Condition requires:
\begin{enumerate}
   \item[1.] $\Phi_{\lambda
\lambda'}^q(X_i')=X_i$ for every  $i>k$;
   \item[2.]$\Phi_{\lambda
\lambda'}^q(X_i')=f_i(X_1,X_2,\cdots,X_k)$ for every $i\leq k$,
where $f_i$ is a multi-variable rational function depending only on the
combinatorics of $\lambda$ and $\lambda'$ in $D(\lambda, \lambda')$ in
the following sense: For any two pairs of ideal triangulation
   $(\lambda, \lambda')$, $(\lambda'', \lambda''')$ for which
there exists a diffeomorphism $\psi: D(\lambda,
\lambda')\rightarrow D(\lambda'', \lambda''')$ sending $\lambda_i$
to $\lambda_j''$ and $\lambda_j'$ to $\lambda_j'''$ for every
$1\leq j\leq k$, then
\begin{align*}
\Phi_{\lambda
\lambda'}^q(X_i') &=f_i(X_1,X_2,\cdots,X_k) \textrm{ and}\\
\Phi_{\lambda''
\lambda'''}^q(X_i''') &=f_i(X_1'',X_2'',\cdots,X_k'')
\end{align*}
for the same
rational function $f_i$.
\end{enumerate}

\section{Diagonal Exchanges}
\label{sect:DiagonalExchange}

The permutation group $\mathcal{S}_n$ acts on the set $\Lambda(S)$ of
isotopy classes of (indexed) ideal triangulations  of
$S$ by reindexing of their components. The set
$\Lambda(S)$ admits another natural operation.

The $i$-th
\textit{diagonal exchange}
$\Delta_i:\Lambda(S) \rightarrow \Lambda(S)$ is defined as
follows. The $i$-th component $\lambda_i$ of the ideal
triangulation $\lambda$ separates two triangle components $T_1$
and $T_2$ of $S_\lambda$. If these two triangles are distinct, the union
$T_1\cup T_2\cup\lambda_i$ is an open square $Q$ with diagonal
$\lambda_i$. Then the ideal triangulation $\Delta_i(\lambda) \in
\Lambda(S)$ is obtained from
$\lambda$ by replacing $\lambda_i$ by the other diagonal $\lambda_i'$ of
the square $Q$. In the remaining case where $T_1=T_2$, then
$\Delta_i(\lambda) =\lambda$ by convention; note that this case occurs
exactly when $\lambda_i$ is the only component of $\lambda$ converging
towards a certain puncture of $S$.

We say that the ideal triangulation $\lambda'=\Delta_i(\lambda)$ is
obtained from $\lambda$ by an \textit{embedded} diagonal exchange if,
with the above notation, the four sides of the square $Q= T_1\cup
T_2\cup\lambda_i$ correspond to distinct components of $\lambda$ (and of
$\lambda'$).

\begin{figure}[h]
\begin{center}
\includegraphics{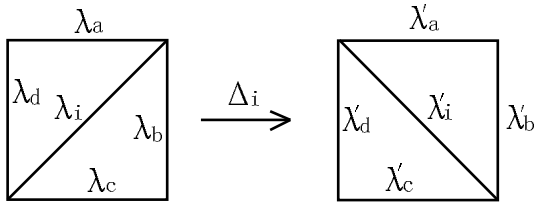}
\caption{A diagonal exchange} \label{DiagonalExchange}
\end{center}
\end{figure}

The coordinate change isomorphisms $\Phi_{\lambda,\lambda'}^q$
constructed by Chekhov and Fock have the property that, if
$\lambda'=\Delta_i(\lambda)$ is obtained from
$\lambda$ by an embedded $i$-th diagonal exchange with indices as shown in
Figure~\ref{DiagonalExchange},
\begin{align*}
\Phi_{\lambda,\lambda'}^q(X_a') & =(1+q X_i)X_a\\
\Phi_{\lambda,\lambda'}^q(X_b') &=(1+q
   X^{-1}_i)^{-1}
X_b \\
   \Phi_{\lambda,\lambda'}^q(X_c')& = (1+q X_i) X_c\\
\Phi_{\lambda,\lambda'}^q(X_d') &=(1+q
   X_i^{-1} )^{-1} X_d,\\
\Phi_{\lambda,\lambda'}^q(X_i') &=  X_i^{-1}
\end{align*}
while $\Phi_{\lambda,\lambda'}^q(X_j')= X_j$ for every $j\not= a$, $b$,
$c$, $d$, $i$.

One needs different formulas for non-embedded diagonal exchanges, which
are given in \cite{Liu04}.

There is an issue that we have somewhat neglected when stating
Theorem~\ref{thm:CoordChangesExist}, which is that it is not completely
trivial that the  coordinate change isomorphisms
$\Phi_{\lambda,\lambda'}^q$ of Chekhov and Fock satisfy the Locality
Condition of Section~\ref{sect:LocalityCondition}. Now is probably a good
time to quickly address this issue.

\begin{lemma}
The Chekhov-Fock coordinate change isomorphisms
$\Phi_{\lambda,\lambda'}^q$ of \cite{Fock97, ChekhovFock00,
ChekhovFock99, Liu04} satisfy the Locality Condition of
Section~\ref{sect:LocalityCondition}.
\end{lemma}
\begin{proof}
Having defined the coordinate change isomorphisms when $\lambda$ and
$\lambda'$ are related by a diagonal exchange or a reindexing,
Chekhov and Fock define $\Phi_{\lambda,\lambda'}^q$ for general $\lambda$,
$\lambda'$ as follows. They connect $\lambda$ to $\lambda'$ by a sequence
$\lambda=\lambda^{(0)}$, $\lambda^{(1)}$, \dots, $\lambda^{(k)}=\lambda'$
such that each $\lambda^{(i)}$ is obtained from $\lambda^{(i-1)}$ by a
diagonal exchange or a reindexing; the existence of such a sequence is
for instance guaranteed by \cite{Penner87, Harer86, Hatcher91}. Then they
define $\Phi_{\lambda,\lambda'}^q$  as
$$
\Phi_{\lambda,\lambda'}^q= \Phi_{\lambda,\lambda^{(1)}}^q\circ
\Phi_{\lambda^{(1)},\lambda^{(2)}}^q \circ \dots
\Phi_{\lambda^{(k-1)},\lambda'}^q
$$
and show that this is independent of the sequence of $\lambda^{(i)}$
connecting $\lambda$ to $\lambda'$. (The possibility of non-embedded
diagonal exchanges in this sequence is neglected by Chekhov
and Fock, and is rigorously dealt with in \cite{Liu04}).

Hatcher \cite{Hatcher91} proves a stronger result, namely that one can
choose the sequence of reindexing and diagonal exchanges in
which all the $\lambda^{(i)}$ are contained in the discrepancy span
$D(\lambda, \lambda')$. If one uses this sequence to define
$\Phi_{\lambda,\lambda'}^q$, it becomes evident that
$\Phi_{\lambda,\lambda'}^q$ satisfies the Locality Condition.
\end{proof}

The main step in the proof of Theorem~\ref{thm:CoordChangesUnique} is the
following, which is its specialization to the case of embedded diagonal
exchanges.

\begin{proposition} [Main step]
\label{MainStep}
Assume that $\chi(S)<-2$. Suppose that we have a
family of algebraic isomorphisms
$$
\Psi_{\lambda,\lambda'}^q:\mathbb{C}(X_1',X_2',\dots,X_n')_{\lambda'}^q
\rightarrow \mathbb{C}(X_1,X_2,\dots,X_n)_{\lambda'}^q
$$ indexed
by pairs of ideal triangulations $\lambda$, $\lambda' \in
   \Lambda(S)$ such that:
\begin{enumerate}
   \item[(a)]
$\Psi_{\lambda,\lambda''}^q=\Psi_{\lambda \lambda'}^q \circ
\Psi_{\lambda',\lambda''}^q$  for any $\lambda$, $\lambda'$ and
$\lambda'' \in \Lambda(S)$ ;
   \item[(b)] if $\lambda'=\sigma
\lambda$ is obtained by re-indexing  $\lambda$ by  $\sigma \in S_n$, then
$\Psi_{\lambda \lambda'}^q(X_i')=X_{\sigma(i)}$ for every $i$;
\item[(c)] the $\Psi_{\lambda,\lambda'}^q$ satisfy the Locality Condition
of Section~\ref{sect:LocalityCondition}.
\end{enumerate}

Then, there exists an invertible element $\xi \in \mathbb C(q)$ and
$\varepsilon = \pm1$, defining for each ideal triangulation $\lambda$
an isomorphism
$$\Theta_{\lambda}:
\mathbb{C}(X_1,X_2,\dots,X_n)_{\lambda}^q \rightarrow
\mathbb{C}(X_1,X_2,\dots,X_n)_{\lambda}^q$$
by the property that $\Theta_{\lambda}(X_k)=\xi X_k^\varepsilon$ for
every
$k$, such that
\begin{align*}
\Theta_{\lambda}^{-1}
\circ \Psi_{\lambda \lambda'}^q \circ
\Theta_{\lambda'}(X_a') & = (1+q X_i)X_a
   \\
\Theta_{\lambda}^{-1}
\circ \Psi_{\lambda \lambda'}^q \circ
\Theta_{\lambda'}(X_b') &=(1+q X^{-1}_i)^{-1} X_b\\
\Theta_{\lambda}^{-1}
\circ \Psi_{\lambda \lambda'}^q \circ
\Theta_{\lambda'}(X_c')& =  (1+q X_i) X_c \\
\Theta_{\lambda}^{-1}
\circ \Psi_{\lambda \lambda'}^q \circ
\Theta_{\lambda'}(X_d') &=(1+q X_i^{-1} )^{-1} X_d \\
\Theta_{\lambda}^{-1}
\circ \Psi_{\lambda \lambda'}^q \circ
\Theta_{\lambda'}(X_i') &= X_i^{-1}
\end{align*}
whenever
$\lambda'=\Delta_i(\lambda)$ is obtained from
$\lambda$ by an embedded $i$-th diagonal exchange with indices as shown
   in Figure~\ref{DiagonalExchange}.
\end{proposition}

The next section is devoted to the proof of
Proposition~\ref{MainStep}.

\section{Proof of Proposition \ref{MainStep}}
\label{ProofMainStep}

We begin with a few algebraic preliminaries which will be regularly used
throughout the proof.

A division algebra $A$ and an algebra isomorphism $\alpha:A\to A$
determine a \emph{skew polynomial algebra} $A[X]_\alpha$ consisting of
all formal polynomials $\sum_{i=1}^n a_i X_i$ in a variable $X$, but
where the multiplication is defined in such a way that $Xa = \alpha(a)
X$. This is a special case of Ore extension \cite{Cohn95, Kassel95,
BrownGoodearl02}. A fundamental property of these skew polynomial
algebras is that they satisfy the Ore condition, so that they can be
enlarged to a fraction division algebra $A(X)_\alpha$. The elements of
$A(X)_\alpha$ are \emph{skew rational fractions}, which can be
expressed as rational fractions in $X$ with coefficients in $A$ but
multiply according to the relation $Xa = \alpha(a)X$.

In particular, the Chekhov-Fock division algebras
$\mathbb{C}(X_1,X_2,\dots,X_n)_{\lambda}^q$ are obtained by iteration
of such constructions.

\begin{lemma}
\label{useful1}
Let $A$ be a division algebra over the field
   $\mathbb{C}(q)$, and consider the skew rational function algebra
$A(X)_\alpha$ associated to an isomorphism $\alpha: A \to A$. If the
skew rational fraction
$f(X)\in A(X)_\alpha$ satisfies
$f(X)=f(q^2X)$, then
$f(X)$ is a constant function.
\end{lemma}
\begin{proof}
As in the commutative case, $f(X)$ admits a unique Laurent series
expansion $f(X)=\sum_{k=n}^{\infty}a_k X^k$ with $n\in \mathbb{Z}$
and $a_k \in A$ for every $k$. Then $f(X)=f(q^2X)$ implies that
$a_k=a_k q^{2k}$ for every $k$.  Hence $a_k=0$ if $k\neq 0$, and
$f(X)$ is a constant function.
\end{proof}

   In the commutative set-up where $\alpha$ is the identity, we will
frequently use the following elementary property, which we
consequently state as a lemma.

\begin{lemma}
\label{useful2}
If $A$ is a division algebra over the field $\mathbb{C}(q)$ and if
$f(X)\in A(X)$ is a (commutative) rational function in the variable $X$,
then
$$
   Xf(Y)X^{-1}=f(XYX^{-1})
$$
in the division algebra $A\{X,Y\}$ consisting of all rational fractions
in the non-commuting variables $X$ and $Y$ (but where $X$ and $Y$
commute with the elements of $A$).
\qed
\end{lemma}

We are now ready to begin the proof of Proposition~\ref{MainStep}.
Suppose that we are given a family of algebraic isomorphisms
$$
\Psi_{\lambda,\lambda'}^q:\mathbb{C}(X_1',X_2',\dots,X_n')_{\lambda'}^q
\rightarrow \mathbb{C}(X_1,X_2,\dots,X_n)_{\lambda}^q
$$ indexed
by pairs of ideal triangulations $\lambda$, $\lambda' \in
   \Lambda(S)$ which satisfy the hypotheses of Proposition~\ref{MainStep}.

Throughout this section, we will assume that $\chi(S)<-2$.

To avoid repetition, whenever we encounter a diagonal exchange, we will
also implicitly assume that the components of the corresponding ideal
triangulations are indexed as in Figure~\ref{DiagonalExchange}.

The next four lemmas use only the
Locality Condition and the fact that the
$\Psi_{\lambda,\lambda'}^q$ are algebra isomorphisms.

\begin{lemma}
\label{perturbation}
There exists rational fractions $f$, $g$ and $h\in\mathbb{C}(q)(X)$ such
that, for every embedded diagonal exchange indexed as in
Figure~\ref{DiagonalExchange},
\begin{align*}
\Psi_{\lambda,\lambda'}^q(X_a') & =f(X_i) X_a\\
\Psi_{\lambda,\lambda'}^q(X_b') &=g(X_i) X_b \\
   \Psi_{\lambda,\lambda'}^q(X_c')& = f(X_i) X_c\\
\Psi_{\lambda,\lambda'}^q(X_d') &=g(X_i) X_d,\\
\Psi_{\lambda,\lambda'}^q(X_i') &= h(X_i).&
\end{align*}
\end{lemma}

\begin{proof}
By the Locality Condition, we know there is a
multi-variable rational function $f_1$ such that
$\Psi_{\lambda,\lambda'}^q (X_a')=f_1(X_a, X_b, X_c, X_d, X_i)$ for
every embedded diagonal exchange. We are going to take advantage of the
fact that $f_1$ is independent of the global properties of $\lambda$ and
$\lambda'$, and depends only on the indexing of the components of
$\lambda$ on the square $D(\lambda, \lambda')$ where the diagonal
exchange takes place.

Because of our hypothesis that $\chi(S)<-2$, we can find by inspection a
pair of ideal triangulations $\lambda$, $\lambda'$ which are deduced from
each other by an embedded diagonal exchange and for which, in addition,
there is a component $\lambda_e$ of $\lambda$ such that the corresponding
generator $X_e$ commutes with $X_a$, $X_b$, $X_c$, $X_i$ and skew-commutes
with $X_d$ in the sense that $X_eX_d = q^2 X_dX_e$. Figure~\ref{not small}
provides examples of such $\lambda$ for the smallest surfaces allowed by
our hypothesis that $\chi(S)<-2$; these actually are the hardest cases
and the reader will easily extend these examples to surfaces with
more punctures.

\begin{figure}[htb]
\begin{center}
\includegraphics{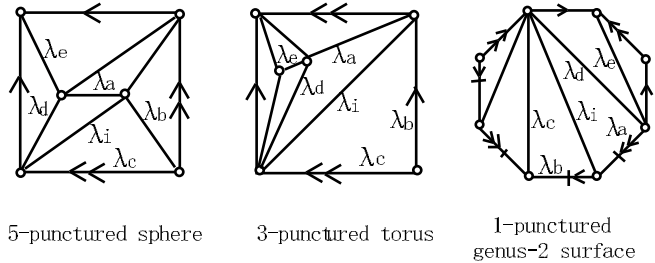}
\caption{Triangulations of sufficiently large surfaces} \label{not
small}
\end{center}
\end{figure}

In this situation, we have $X_e'X_a'{X_e'}^{-1}=X_a'$. Therefore
$$\Psi_{\lambda,\lambda'}^q
(X_e')\Psi_{\lambda,\lambda'}^q (X_a')\Psi_{\lambda,\lambda'}^q
(X_e')^{-1}=\Psi_{\lambda,\lambda'}^q (X_a')$$ since
$\Phi_{\lambda,\lambda'}^q$ is  an algebra isomorphism. Because
the components $\lambda_e$ and $\lambda_e'$ are not in
$D(\lambda,\lambda')$, the Locality Condition implies that
$\Psi_{\lambda,\lambda'}^q (X_e')=X_e$. It follows that
$X_ef_1(X_a,X_b,X_c,X_d,X_i)X_e^{-1}=f_1(X_a,X_b,X_c,X_d,X_i)$.
Applying Lemma~\ref{useful2}, we conclude that
$f_1(X_a,X_b,X_c,q^2X_d,X_i)=f_1(X_a,X_b,X_c,X_d,X_i)$.

Lemma~\ref{useful1} then shows that
$f_1$ must be independent of the variable $X_d$.

Similarly, one can find another pair of ideal triangulations $\lambda$
and $\lambda'$, obtained from each other by a diagonal exchange, for
which there is now a generator $X_e$ which skew-commutes with
$X_c$ and commutes with the generators $X_a$, $X_b$, $X_d$ associated to
the other sides of the square
$D(\lambda, \lambda')$. The same argument as above then shows that $f_1$
is independent of $X_c$.

Another application of the same argument shows that $f_1$ is independent
of $X_b$.

Therefore,  $\Phi(X_a')=f_1(X_a,X_b,X_c,X_d,X_i)=f_2(X_a,X_i)$ for some
rational fraction in two variables.

Now, we again find two ideal triangulations $\lambda$
and $\lambda'$, obtained from each other by a diagonal exchange, for
which there is now a generator $X_e$ which commutes
$X_b$, $X_c$, $X_d$ and for which $X_eX_a = q^2 X_aX_e$. The argument is
here slightly different.
From
$X_e'X_a'{X_e'}^{-1}=q^2X_a'$, we deduce
$f_2(q^2X_a,X_i)=q^2f_2(X_a,X_i)$. This is equivalent to say
$f_2(q^{-2}X_a,X_i)(q^2 X_a)^{-1} =f_2(X_a,X_i)X_a^{-1}$, which forces
$f_2(X_a,X_i)X_a^{-1}$ to be independent of  variable $X_a$ by
Lemma~\ref{useful1}. In other words,
$\Psi_{\lambda,\lambda'}^q (X_a')=f_2(X_a,X_i)= f(X_i) X_a$ for some
rational function $f \in
\mathbb C(q)$, as required.

Combining the Locality Condition with the symmetries of the
square, $\Psi_{\lambda,\lambda'}^q (X_b')= g(X_i) X_b$,
$\Psi_{\lambda,\lambda'}^q (X_c')= f(X_i) X_c$, and
$\Psi_{\lambda,\lambda'}^q (X_d')= g(X_i) X_d$ for another
rational function $g \in \mathbb C(q)$.

Finally, we know by the Locality Condition that
$\Psi_{\lambda,\lambda'}^q(X_i')=h_1(X_a,X_b,X_c,X_d,X_i)$ for some
multi-variable rational function $h_1$. Finding again a pair
$\lambda$,
$\lambda'$ where there is a generator $X_e$ which commutes with
$X_b$, $X_c$, $X_d$ ,$X_i$ but skew-commutes with $X_a$, the same method
as above shows that $\Psi_{\lambda,\lambda'}^q(X_i')$ is independent of
$X_a$. Similarly
$\Psi_{\lambda,\lambda'}^q(X_i')$ is also independent of $X_b$, $X_c$ and
$X_d$, hence
$\Psi_{\lambda,\lambda'}^q(X')=h(X)$ for some rational function $h \in
\mathbb C(q)$.
\end{proof}

\begin{lemma}
\label{Function h}
The rational function $h$ of Lemma~\ref{perturbation} is
of the form $h(X)= \eta X^{-1}$ for some $\eta
\in \mathbb{C}(q) -\{0\}$.
\end{lemma}

\begin{proof}
Applying $\Psi_{\lambda,\lambda'}^q$ to both sides of the relation
$X_i' X_a'=q^2 X_a' X_i'$, we obtain that $ h(X_i)f(X_i) X_a=q^2
f(X_i) X_a h(X_i)$. The terms $f(X_i)$ cancel out. Using the
relation $X_a X X_a^{-1}=q^2 X$ and Lemma~\ref{useful2}, we get
$h(q^{-2}X_i)=q^2 h(X_i)$. Rewriting this as
$(q^{-2}X_i)\,h(q^{-2}X_i)=X_i h(X_i)$, we see that $X_i
h(X_i)=\eta$ for some constant $\eta \in \mathbb{C}(q)$ by Lemma
\ref{useful1}. Because the subalgebra of
$\mathbb{C}(X_1,X_2,\dots,X_n)_{\lambda}^q$ consisting of all
rational functions in the variable $X_i$ is abstractly isomorphic
to $\mathbb C(q)(X)$, we conclude that $h(X)=\eta X^{-1}$ in
$\mathbb C(q)(X)$.
\end{proof}

The next lemma uses the fact that the diagonal exchange map $\Delta_i
:\Lambda(S) \to \Lambda(S)$ is an involution.

\begin{lemma}
\label{f_and_g}
The rational functions $f$, $g\in \mathbb C(q)(X)$ of
Lemma~\ref{perturbation} are such that $g(X) = f(\eta X^{-1})^{-1}$,
where
$\eta\in  \mathbb C(q)$ is the constant of Lemma~\ref{Function h}.
\end{lemma}
\begin{proof}
Because $\lambda = \Delta_i(\lambda') = \Delta_i^2(\lambda)$,
\begin{align*}
X_a' &= \Psi_{\lambda',\lambda}^q \circ \Psi_{\lambda,\lambda'}^q (X_a')
= \Psi_{\lambda',\lambda}^q(f(X_i)X_a)
= f(\eta {X_i'}^{-1}) g(X_i')X_a'
\end{align*}
since $\Psi_{\lambda',\lambda}^q(X_i) = \eta {X_i'}^{-1}$ and
$\Psi_{\lambda',\lambda}^q(X_a)= g(X_i')X_a'$. The result follows.
\end{proof}

\begin{lemma}
\label{f}  The rational function $f$ in  Lemma \ref{perturbation}
   satisfies $f(X)=\xi qX f(q^{-2} \xi^{-2}X^{-1})$ for
some $\xi \in \mathbb C(q)$ with $\xi^2 = \eta^{-1}$.
\end{lemma}

\begin{proof}
The product $X_a'X_b'X_c'X_d'X_i'$  commutes with each of
$X_a'$, $X_b'$, $X_c'$, $X_d'$ and $X_i'$. Since
$\Psi_{\lambda',\lambda}^q$ is an algebra isomorphism,
$\Psi_{\lambda',\lambda}^q(X_a'X_b'X_c'X_d'X')$ therefore commutes with
   $X_a$, $X_b$, $X_c$, $X_d$ and $X_i$. Using Lemma~\ref{useful2} and the
skew-commutativity relations, this element is equal to
\begin{align*}
\Psi_{\lambda',\lambda}^q(X_a'X_b'X_c'X_d'X_i') &
=f(X_i)X_a \cdot f(\eta X_i^{-1})^{-1}X_b
\cdot f(X_i) X_c \cdot f(\eta X_i^{-1})^{-1}X_d \cdot \eta X_i^{-1}
\\ &= \eta X_d X_a X_b X_c X_i \left[f(X_i)(qX_i f(q^{-2}\eta
X_i^{-1}))^{-1}\right]^{2}
\end{align*}
Since $X_d X_a X_b X_c X_i$ commutes with $X_a$, $X_b$, $X_c$,
$X_d$ and $X_i$, it follows that \linebreak $\left[ f(X_i)(qX_i
f(q^{-2} \eta X_i^{-1}))^{-1}\right]^2$ also
   commutes with $X_a$.  An application of Lemmas~\ref{useful2}
and \ref{useful1} then shows that this rational function of $X_i$ is
constant. As a consequence,
$f(X_i)(qX_i f(q^{-2} \eta
X_i^{-1}))^{-1}= \xi$ for some
$\xi\in \mathbb{C}(q)$.

Therefore, $f(X)(qX f(q^{-2} \eta
X^{-1}))^{-1}= \xi$ in $\mathbb C (q)(X)$.
Substituting $q^{-2}\eta X^{-1}$ for $X$ in this equation and combining
with the original relation, we obtain that $\xi^2 = \eta^{-1}$.

This proves the relation $f(X)=\xi qX f(q^{-2}
\xi^{-2}X^{-1})$.
\end{proof}

Without loss of generality, we can replace all isomorphisms
$$
\Psi_{\lambda
\lambda'}^q:\mathbb{C}(X_1',X_2',\dots,X_n')_{\lambda'}^q
\rightarrow \mathbb{C}(X_1,X_2,\dots,X_n)_{\lambda}^q
$$
by
$\Theta_{\lambda}^{-1}
\circ \Psi_{\lambda \lambda'}^q \circ
\Theta_{\lambda'}$ where, for each ideal triangulation $\lambda$,
$$\Theta_{\lambda}:
\mathbb{C}(X_1,X_2,\dots,X_n)_{\lambda}^q \rightarrow
\mathbb{C}(X_1,X_2,\dots,X_n)_{\lambda}^q$$
is the isomorphism
defined by $\Theta_{\lambda}(X_k)=\xi X_k$ for every $k$. This replaces
the rational function $f(X)$ by $f'(X) = f(\xi^{-1}X)$ which satisfies
the relation $f'(X) = qX f'(q^{-2}X^{-1})$, the function $g(X) =
f(\xi^{-2} X^{-1})^{-1}$ by $g'(X) = f'(X^{-1})^{-1}$, and $h(X) =
\xi^{-2} X^{-1}$ by
$h'(X) = X^{-1}$. Therefore, we can henceforth assume that $\xi=1$.

\begin{figure}[h]
\begin{center}
\includegraphics{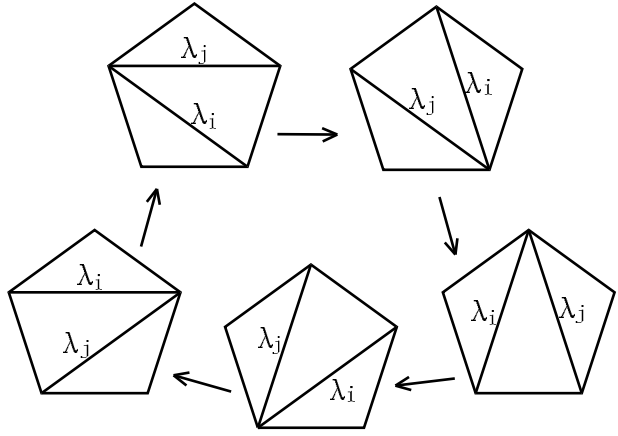}
\caption{The Pentagon Relation}
\label{Pentagon Move}
\end{center}
\end{figure}

So far, we have only used the Locality Condition and the
skew commutativity properties of the generators of the Chekhov-Fock
algebra. We now take advantage a subtler property, namely the Pentagon
Relation illustrated in Figure~\ref{Pentagon Move}. Because of our
assumption that
$\chi(S)<-2$, the surface $S$ admits an ideal triangulation $\lambda$
for which the closure of three components of $S-\lambda$ forms an embedded
pentagon. We then have a sequence of 5 diagonal exchanges which returns
to the original ideal triangulation. (Actually, the indexing is modified
by a transposition of the  two diagonals of the pentagon, but this will
be irrelevant here).

For each ideal triangulation $\lambda$ occurring in the pentagon move, it
is convenient to index the components $\lambda_i$ and $\lambda_j$ of
$\lambda$ that are diagonals of the pentagon as indicated on
Figure~\ref{Pentagon Move}. This convention has the advantage that the
corresponding generators $X_i$ and $X_j$ of the Chekhov-Fock
algebra of $\lambda$ satisfy the same relation
$X_jX_i = q^2 X_iX_j$. Also, as one moves from one such ideal
triangulation $\lambda$ to the next one $\lambda'$,
$\Psi_{\lambda,\lambda'}^q (X_i')= f(X_i)X_j$ and
$\Psi_{\lambda,\lambda'}^q (X_j')=X_i^{-1}$.

This leads us to introduce the algebra $\mathcal W$ defined by generators
$U^{\pm1}$, $V^{\pm1}$ and by the relation $VU=q^2UV$, and the algebra
isomorphism $\Psi: \mathcal W \to \mathcal W$ defined by $\Psi(U) = f(U)V$
and
$\Psi(V) = U^{-1}$. From the above observation, a consequence of the
Pentagon Relation is that the composition $\Psi^5$ is the identity.

\begin{lemma}
\label{order5} The isomorphism $\Psi: \mathcal W \to \mathcal W$  defined
by $\Psi(U) = f(U)V$ and $\Psi(V) = U^{-1}$ has period $5$ if and only if
$f(f(U)V)=Vf(U^{-1}f(V^{-1}))U$.
\end{lemma}

\begin{proof} Write $U_k=\Psi^k(U)$ and $V_k = \Psi^k(V)$. Because
$V_{k+1}=U_k^{-1}$,  $\Psi$ is periodic with period 5 if and only if
$U_{k+5}=U_k$ for some $k$ or, equivalently, or for all $k$.

Noting that $U_{k+1} = f(U_k)V_k = f(U_k)U_{k-1}^{-1}$ and
$U_{k-1}=U_{k+1}^{-1}f(U_k)$, an immediate computation gives
$U_2=f(f(U)V)U^{-1}$ and $U_{-3} = Vf(U^{-1}f(V^{-1}))$. The result
follows.
\end{proof}

\begin{lemma}
\label{f when xi=1}
If a rational function $f$ is such that $f(X)=qXf(q^{-2}X^{-1})$ in
$\mathbb{C}(q)(X)$ and $ f(f(U)V)=V f(U^{-1}f(V^{-1}))U$ in the
algebra $\mathcal W$,
then $f(X)=1+qX$ or $f(X)=(1+q^{-1}X^{-1})^{-1}$.
\end{lemma}

\begin{proof}
Let us transform the relation $f(f(U)V)=Vf(U^{-1}f(V^{-1}))U$.
\begin{align*}
f(f(U)V)&=Vf(U^{-1}f(V^{-1}))U&\\
&=V f\left(q U^{-1} V^{-1}f(q^{-2}V)\right)U&\\
&=V q^2 U^{-1} V^{-1}f(q^{-2}V) f\left(q^{-2} q^{-1}f(q^{-2}V)^{-1}
VU\right)U&\\
&=\left( U^{-1}f(q^{-2}V)U\right) \left(
U^{-1}f(q^{-3}f(q^{-2}V)^{-1} VU)U\right)&\\
&=f(V)f\left( q^{-1}f(V)^{-1}VU\right)
\end{align*}

Therefore the equation $ f(f(U)V)=V f(U^{-1}f(V^{-1}))U$
is equivalent to the relation $f(f(U)V)=f(V)f( q^{-1}f(V)^{-1}VU)$.

Consider the Laurent expansion $f(X)=X^s (a_0+a_1 X+a_2X^2+\dots)$ of
$f$, with $a_0\not =0$. We will distinguish cases.

\medskip
\noindent Case 1: $s<0$.

   If we expand $q^{-1}f(V)^{-1}VU$ as a Laurent series in $V$ with
coefficients in $\mathbb C(q)(U)$, its lowest degree term in $V$ has
degree
$1-s\geq 2$. We can therefore use again the Laurent series expansion of
$f(X)$ to expand both sides of the equation $f(f(U)V)=f(V)f(
q^{-1}f(V)^{-1}VU)$ as a Laurent series in $V$ with
coefficients in $\mathbb C(q)(U)$.

In this expansion, the lowest degree term in $V$ has degree $s+s(1-s)$
for the right hand side, and $s$ of the left hand side. Since
$s=s+s(1-s)$ has no negative solution, it follows that this case
actually cannot occur.

\medskip\noindent
   Case $2$: $s=0$.

In this case the lowest degree term of $q^{-1}f(V)^{-1}VU$ has degree 1,
so we can expand again $f( q^{-1}f(V)^{-1}VU)$ term by term.

   Comparing the constant terms (with respect to $V$) on
   both sides of $f(f(U)V)=f(V)f( q^{-1}f(V)^{-1}VU)$, we get
   $a_0=a_0^2$. So $a_0=1$.

   The next term in the expansion has degree $t=\min\{T: T>0, a_T\ne
0\}$. The coefficient of
   $V^{t}$ on the left hand side of the equation comes from
$a_{t}(f(U)V)^t$,
   and is equal to $a_tf(U)f(q^{2}U)\dots f(q^{2(t-1)}U)$. On the
   right hand side, remembering that $a_0=1$, the coefficient of $V^t$
comes from $a_t V^t + a_t (q^{-1}VU)^t$  and is equal to $a_t +
a_t q^{t^2}U^t$. This gives
$$f(U)f(q^{2}U)\dots f(q^{2(t-1)}U) = 1 +  q^{t^2}U^t.$$
Expanding the left hand side as a Laurent series in $U$ and looking at
the coefficient of $U^t$, we conclude that
$t=1$ and $f(X) = 1+ qU$.

\medskip\noindent
Case 3: $s>0$.

In this case, the Laurent expansion of $q^{-1}f(V)^{-1}VU$ in $V$ begins
with a term of non-positive degree, so that we cannot expand $f(
q^{-1}f(V)^{-1}VU)$ term by term any more. We will use an internal
symmetry of the problem.

Set $\widehat f
(X)=f(q^{-2}X^{-1})^{-1}=qX f(X)^{-1}$. We claim that this rational
function
$\widehat f
(X)$ also satisfies the hypotheses of the Lemma. Let us check this.
First,
$$qX\widehat f
(q^{-2}X^{-1})=qXf(X)^{-1}=f(q^{-2}X^{-1})^{-1}=\widehat f
(X).$$
Then
\begin{align*}
\widehat f
(\widehat f
(U)V)& =f(q^{-2}V^{-1}\widehat f
(U)^{-1})^{-1}=
f(q^{-2}V^{-1}f(q^{-2}U^{-1}))^{-1}\\
&= f(f(U^{-1})q^{-2}V^{-1})^{-1},
\end{align*}
   while
\begin{align*}
\widehat f
(V)\widehat f
(q^{-1}\widehat f
(V)^{-1}VU)&=f(q^{-2}V^{-1})^{-1}f(q^{-2}qU^{-1}V^{-1}\widehat f
(V))^{-1}\\
&=f(q^{-2}V^{-1})^{-1}f\left(q^{-1}U^{-1}V^{-1}
f(q^{-2}V^{-1})^{-1}\right)^{-1}\\
&=f\left(q^{-1}f(q^{-2}V^{-1})^{-1}q^{-2}V^{-1}U^{-1}\right)
^{-1}f(q^{-2}V^{-1})^{-1}.
\end{align*}
Therefore we need to check that
$$f(f(U^{-1})q^{-2}V^{-1})= f(q^{-2}V^{-1})
f\left(q^{-1}f(q^{-2}V^{-1})^{-1}q^{-2}V^{-1}U^{-1}\right).$$ But this
property holds by applying  to the identity
$f(f(U)V)=f(V)f( q^{-1}f(V)^{-1}VU)$ the homomorphism $\Theta:
\mathcal W\rightarrow \mathcal W$ defined by
$\Theta(U)=U^{-1}$ and $\Theta(V)=q^{-2}V^{-1}$.

Therefore, the rational function $\widehat f
(X) \in \mathbb C(q)(X)$
satisfies the hypotheses of the Lemma. The Laurent series
expansion of $\widehat f
(X)= qXf(X)^{-1}$ begins with a term of degree
$1-s\geq 0$. It follows from Cases 1 and 2 that $\widehat f
(X)= 1+qX$.
Therefore, $f(X) = qX\widehat f
(X)^{-1} = (1+q^{-1}X^{-1})^{-1}$.
   \end{proof}

Note that, when $f(X) = 1+qX$, the function $g(X)$ of
  Lemmas~\ref{perturbation} and
\ref{f_and_g} is such that $g(X) = f(X^{-1})^{-1}= (1+qX^{-1})^{-1}$.
When $f(X) =
(1+q^{-1}X)^{-1}$, then
$g(X)=1+q^{-1}X$.

At this point we have shown that there exists  $\xi \in \mathbb
C-\{0\}$ such that, if isomorphisms
$$\Theta_{\lambda}:
\mathbb{C}(X_1,X_2,\dots,X_n)_{\lambda}^q \rightarrow
\mathbb{C}(X_1,X_2,\dots,X_n)_{\lambda}^q$$ are defined by
$\Theta_{\lambda}(X_k)=\xi X_k$ for every $k$, either we have
\begin{align*}
\Theta_{\lambda}^{-1}
\circ \Psi_{\lambda \lambda'}^q \circ
\Theta_{\lambda'}(X_a') & = (1+q X_i)X_a
   \\
\Theta_{\lambda}^{-1}
\circ \Psi_{\lambda \lambda'}^q \circ
\Theta_{\lambda'}(X_b') &=(1+q X^{-1}_i)^{-1} X_b\\
\Theta_{\lambda}^{-1}
\circ \Psi_{\lambda \lambda'}^q \circ
\Theta_{\lambda'}(X_c')& =  (1+q X_i) X_c \\
\Theta_{\lambda}^{-1}
\circ \Psi_{\lambda \lambda'}^q \circ
\Theta_{\lambda'}(X_d') &=(1+q X_i^{-1} )^{-1} X_d \\
\Theta_{\lambda}^{-1}
\circ \Psi_{\lambda \lambda'}^q \circ
\Theta_{\lambda'}(X_i') &= X_i^{-1}
\end{align*}
whenever
$\lambda'$ is obtained from
$\lambda$ by an embedded diagonal exchange (with the edge indicing
conventions of  Figure~\ref{DiagonalExchange}), or we have
\begin{align*}
\Theta_{\lambda}^{-1}
\circ \Psi_{\lambda \lambda'}^q \circ
\Theta_{\lambda'}(X_a') & = (1+q^{-1} X_i^{-1})^{-1}X_a
   \\
\Theta_{\lambda}^{-1}
\circ \Psi_{\lambda \lambda'}^q \circ
\Theta_{\lambda'}(X_b') &=(1+q^{-1} X_i) X_b\\
\Theta_{\lambda}^{-1}
\circ \Psi_{\lambda \lambda'}^q \circ
\Theta_{\lambda'}(X_c')& =  (1+q^{-1} X_i^{-1})^{-1} X_c \\
\Theta_{\lambda}^{-1}
\circ \Psi_{\lambda \lambda'}^q \circ
\Theta_{\lambda'}(X_d') &=(1+q^{-1} X_i ) X_d \\
\Theta_{\lambda}^{-1}
\circ \Psi_{\lambda \lambda'}^q \circ
\Theta_{\lambda'}(X_i') &= X_i^{-1}
\end{align*}
again whenever
$\lambda'$ is obtained from
$\lambda$ by an embedded  diagonal exchange.

In the first case, we have reached the conclusion desired for
Proposition~\ref{MainStep}.

In the second case, we conjugate by an additional family of isomorphisms
   $$\Theta_{\lambda}':
\mathbb{C}(X_1,X_2,\dots,X_n)_{\lambda}^q \rightarrow
\mathbb{C}(X_1,X_2,\dots,X_n)_{\lambda}^q$$
defined by $\Theta_{\lambda}(X_k)= X_k^{-1}$ for every $k$. Then
\begin{align*}
(\Theta_{\lambda}\circ
\Theta_{\lambda}')^{-1}
\circ \Psi_{\lambda \lambda'}^q \circ
(\Theta_{\lambda'}\circ
\Theta_{\lambda'}')(X_a') & = X_a (1+q^{-1} X_i) = (1+q X_i)X_a
   \\
(\Theta_{\lambda}\circ
\Theta_{\lambda}')^{-1}
\circ \Psi_{\lambda \lambda'}^q \circ
(\Theta_{\lambda'}\circ
\Theta_{\lambda'}')(X_b') &=X_b (1+q^{-1} X^{-1}_i)^{-1}  = (1+q
X^{-1}_i)^{-1} X_b\\
(\Theta_{\lambda}\circ
\Theta_{\lambda}')^{-1}
\circ \Psi_{\lambda \lambda'}^q \circ
(\Theta_{\lambda'}\circ
\Theta_{\lambda'}')(X_c')& =  X_c(1+q^{-1} X_i)  =  (1+q X_i) X_c\\
(\Theta_{\lambda}\circ
\Theta_{\lambda}')^{-1}
\circ \Psi_{\lambda \lambda'}^q \circ
(\Theta_{\lambda'}\circ
\Theta_{\lambda'}')(X_d') &=X_d (1+q^{-1} X_i^{-1} )^{-1}  =(1+q
X_i^{-1} )^{-1} X_d \\ (\Theta_{\lambda}\circ
\Theta_{\lambda}')^{-1}
\circ \Psi_{\lambda \lambda'}^q \circ
(\Theta_{\lambda'}\circ
\Theta_{\lambda'}')(X_i') &= X_i^{-1}.
\end{align*}
We again have reached the conclusions of Proposition~\ref{MainStep}.

This completes the proof of this statement. \qed

\section{Proof of the main Theorem 2}

Let $
\Phi_{\lambda
\lambda'}^q:\mathbb{C}(X_1',X_2',\dots,X_n')_{\lambda'}^q
\rightarrow \mathbb{C}(X_1,X_2,\dots,X_n)_{\lambda}^q
$
be the Chekhov-Fock family of isomorphisms provided by Theorem~1, and
let $
\Psi_{\lambda
\lambda'}^q:\mathbb{C}(X_1',X_2',\dots,X_n')_{\lambda'}^q
\rightarrow \mathbb{C}(X_1,X_2,\dots,X_n)_{\lambda}^q
$ be another family of isomorphisms satisfying the same properties.

Proposition~\ref{MainStep} shows that, after conjugating the $\Psi_{\lambda
\lambda'}^q$ with appropriate rescaling/inversion isomorphisms
$
\Theta_{\lambda}:\mathbb{C}(X_1,X_2,\dots,X_n)_{\lambda}^q
\rightarrow \mathbb{C}(X_1,X_2,\dots,X_n)_{\lambda}^q
$ defined by $\Theta_{\lambda}(X_k) = \xi X_k^{\pm1}$, we can arrange
that $\Psi_{\lambda
\lambda'}^q=\Phi_{\lambda
\lambda'}^q$ whenever $\lambda$ and $\lambda'$ differ only by an
embedded diagonal exchange.

Because of our hypothesis that $\chi(S)<-2$, there are only two
possible types of non-embedded diagonal exchanges:
\begin{enumerate}
\item [(1)] one type where, in the square where the diagonal exchange
takes, two opposite sides are identified in $S$;
\item [(2)] another type where  two adjacent
sides of this square are identified; this is possible
only when $S$ has at least two punctures.
\end{enumerate}
Explicit formulas for the
   $\Phi_{\lambda \lambda'}^q$ associated to non-embedded diagonal
exchanges are given in \cite{Liu04}.

By inspection on the surface $S$ and using the fact that $\chi(S)<-2$,
one can construct an ideal triangulation of $S$ where three faces form
a pentagon $P$, in such a way that the only two sides of the pentagon
$P$ that are identified in $S$ are two non-adjacent sides. This gives
rise to a sequence of diagonal exchanges
$$\lambda_{(0)} \to \lambda_{(1)} \to \lambda_{(2)} \to
\lambda_{(3)}
\to \lambda_{(4)} \to \lambda_{(5)} = \lambda_{(0)}, $$
all taking place in the pentagon $P$ as in the pentagon relation, where
$\lambda_{(0)}
\to
\lambda_{(1)}$ is a non-embedded diagonal exchange of Type (1), and
where all other diagonal exchanges $\lambda_{(i)} \to \lambda_{(i+1)}$
are embedded.

  From the properties that $\Phi_{\lambda \lambda''}^q= \Phi_{\lambda
\lambda'}^q \circ
\Phi_{\lambda' \lambda''}^q$ and $\Psi_{\lambda
\lambda''}^q=\Psi_{\lambda \lambda'}^q \circ
\Psi_{\lambda' \lambda''}^q$,  we conclude that
$$
\Phi_{\lambda_{(0)} \lambda_{(1)}}^q =
\left(
\Phi_{\lambda_{(1)} \lambda_{(2)}}^q \circ
\Phi_{\lambda_{(2)} \lambda_{(3)}}^q \circ
\Phi_{\lambda_{(3)} \lambda_{(4)}}^q \circ
\Phi_{\lambda_{(4)} \lambda_{(5)}}^q
\right)^{-1}
$$
and
$$
\Psi_{\lambda_{(0)} \lambda_{(1)}}^q =
\left(
\Psi_{\lambda_{(1)} \lambda_{(2)}}^q \circ
\Psi_{\lambda_{(2)} \lambda_{(3)}}^q \circ
\Psi_{\lambda_{(3)} \lambda_{(4)}}^q \circ
\Psi_{\lambda_{(4)} \lambda_{(5)}}^q
\right)^{-1}.
$$
Since $\Phi_{\lambda\lambda'}^q$ and $\Psi_{\lambda\lambda'}^q$
coincide on embedded diagonal exchanges, we conclude that
$\Phi_{\lambda_{(0)} \lambda_{(1)}}^q =
\Psi_{\lambda_{(0)} \lambda_{(1)}}^q$. By the Locality Condition, it
follows that $\Phi_{\lambda\lambda'}^q$ and $\Psi_{\lambda\lambda'}^q$
coincide on all non-embedded diagonal exchanges of Type~(1).

Similarly, when $S$ has at least two punctures, one can find an ideal
triangulations with a pentagon whose only side identifications are
between two adjacent sides. The same argument as above then shows that
$\Phi_{\lambda\lambda'}^q$ and $\Psi_{\lambda\lambda'}^q$
coincide on all non-embedded diagonal exchanges of Type~(2).

Therefore, $\Phi_{\lambda\lambda'}^q$ and $\Psi_{\lambda\lambda'}^q$
coincide on all  diagonal exchanges, embedded or non-embedded.

Finally, any two ideal triangulations $\lambda$ and $\lambda'$ can be
joined by a sequence of diagonal exchanges
$$\lambda = \lambda_{(0)} \to \lambda_{(1)} \to \dots \to \lambda_{(n)}
=\lambda'. $$
See for instance \cite{Penner87, Hatcher91}. This decomposes
$\Phi_{\lambda\lambda'}^q$ and $\Psi_{\lambda\lambda'}^q$ as a product
of
$\Phi_{\lambda_{(i)}\lambda_{(i+1)}}^q=
\Psi_{\lambda_{(i)}\lambda_{(i+1)}}^q$ associated to diagonal exchanges,
thereby showing that
$\Phi_{\lambda\lambda'}^q =\Psi_{\lambda\lambda'}^q$ for any ideal
triangulations $\lambda$ and $\lambda'$.

This completes the proof of Theorem~2.

\newcommand{\journalname}[1]{\textrm{#1}}
\newcommand{\booktitle}[1]{\textrm{#1}}

\end{document}